%

\magnification=\magstep1
\input amstex
\documentstyle{amsppt}
\topmatter
\title Shooting a club with finite conditions \endtitle
\author Jind\v rich Zapletal \endauthor
\affil The Pennsylvania State University \endaffil
\address Department of Mathematics, The Pennsylvania State 
University, University Park, PA 16802 \endaddress
\email zapletal\@ math.psu.edu \endemail
\abstract We study cohabitation of the poset shooting a club through 
a given stationary subset of $\omega _1$ with finite conditions with 
other forcings. \endabstract

\endtopmatter
\document
\proclaim
{Definition 1} \roster
\item If $I\subset \omega _1$ is a countable interval of ordinals and 
$S\subset \omega _1$ is countable, we define $P_{S,I}=\{ p:p$ is a 
finite function from $I$ to $S$ such that there exists $f:I \to S\cap \alpha$
 increasing continuous, $ rng(f)$ unbounded in $S,$ $p\subset f\} .$
\item Let $S\subset \omega _1$ be stationary, $I=\omega _1 \setminus \alpha$ 
for some $\alpha <\omega _1.$ Then $P_{S,I}=\{ p:p$ is a finite function 
from $I$ to $S$ such that $\exists \alpha <\beta <\omega _1$ $\exists f:
\beta \setminus \alpha \to S$ continuous increasing and $p\subset f\} .$
 The order is by inverse inclusion. $P_S=P_{S,\omega _1}.$
\endroster
\endproclaim
We will be interested in $P_S$ for various $S\subset \omega _1$ 
stationary, ``shooting a club through $S$ with finite conditions". 
\proclaim
{Lemma 1} If $p=\{ \langle \alpha ,\alpha \rangle \} \in P_S$ 
then $P_S\restriction p=P_{S\cap \alpha ,\alpha }\times 
P_{S\setminus \alpha ,\omega _1\setminus \alpha }.$
\endproclaim
 
\proclaim
{Corollary 1} $P_S$ is $S$-proper.\cite {S}
\endproclaim
\demo
{Proof} Let $p\in P_S,$ $N\prec H_\theta$ countable with $S,p\in N$ and 
$\alpha =N\cap \omega _1\in S.$ Then $q=p\cup \{ \langle \alpha ,\alpha 
\rangle \} \in P$ is a master condition for $N$ as in \cite {B} .
\enddemo
\proclaim
{Corollary 2} $P_S$ is homogeneous.
\endproclaim
\demo
{Proof} If $p,q\in P_S$ find $N\prec H_\theta$ countable such that 
$p,q,S\in N$ and $\alpha =N\cap \omega _1\in S.$ Then $p_1=p\cup 
\{ \langle \alpha ,\alpha \rangle \},$ $q_1=q\cup \{ \langle \alpha ,
\alpha \rangle \}$ are both in $P_S$ and due to the Lemma can be viewed 
as elements of $P_{S\cap \alpha ,\alpha }\times P_{S\setminus \alpha ,
\omega _1\setminus \alpha }$ with support only the first coordinate. 
But $P_{S\cap \alpha ,\alpha }$ is a countable notion of forcing and 
so is homogeneous. Now it is easy to devise an automorphism of $P_S$ 
sending $p_1$ under $q_1,$ proving homogeneity.
\enddemo
Fix $G\subset P_S$ generic.

\proclaim
{Corollary 3} $r\in \omega ^\omega \cap V[G]$ iff $\exists \alpha <
\omega _1$ $r\in V[G\restriction \alpha ].$
\endproclaim

\proclaim
{Corollary 4} If $V\subset W,$ $\omega _1 ^V=\omega _1 ^W$ and $G\in W$ 
then $G\subset P_S$ is $V$-generic iff \roster
\item $G:\omega _1\to S$ is increasing and continuous
\item $\forall \alpha <\omega _1$ $G(\alpha )=\alpha$ implies 
$G\restriction \alpha$ is $P_{S\cap \alpha ,\alpha }$-generic.
\endroster
\endproclaim
\demo
{Proof} Let $G$ has the above properties and $A\subset P_S$ is a 
maximal antichain in $V.$ Choose $N\prec H_\theta$ in $W$ countable 
containing $G,A.$ Then $A\cap N$ is a maximal antichain in $N\cap P_S$ 
and so it is met by (2). (1) actually follows from (2).
\enddemo

\proclaim
{Lemma 2} If $T$ is a tree of height $\omega _1$ and $P_S\ni p\Vdash$
``$\dot b$ is a branch through $T"$ then $p\Vdash \dot b\in V.$
\endproclaim

\demo
{Proof} Let $\theta$ be large enough regular cardinal and $M\prec H_\theta$
be a countable submodel containing $S,p,T,\dot b,$ $\alpha =M\cap \omega _1
\in S.$   Set $p_1=p\cup \{ \langle \alpha ,\alpha \rangle \} .$ Find 
$p_2\leq p_1$ deciding $\dot b\cap$ the $\beta ^{th}$ level of $T.$ 
Let $q=p_2\cap M.$ For any $r_1,r_2\in M,$
if $r_1,r_2\leq q$ then both of them are comparable with $p_2$ and therefore
any elements of $T\cap M$ forced by $r_1$ or $r_2$ into $\dot b$ must
be linearly ordered in $T.$ By elementarity, $c=\{ x\in T:\exists r\leq q$
$r\Vdash x\in \dot b\}$ is a branch of $T$ and consequently $q\Vdash \dot b=c.$
Since $q\leq p$ and $c\in V$ we are done.
\enddemo

\proclaim
{Lemma 3} If $Q$ is c.c.c. then $P_S\Vdash$``$Q$ is c.c.c.".
\endproclaim

\demo
{Proof} Let $p\Vdash$``$\langle \dot q_\alpha :\alpha <\omega _1\rangle$ 
is an antichain in $Q".$ Fix a bijection $g:P_S\to \omega _1$ and find 
$\langle M_\alpha :\alpha <\omega _1\rangle ,$ a continuous increasing 
chain of submodels of $H_\theta$ with $S,g,p,\langle \dot q_\alpha :
\alpha <\omega _1\rangle \in M_0.$ Set $\beta _\alpha =M_\alpha \cap 
\omega _1$ and find $\langle \dot p_\alpha :\alpha <\omega _1\rangle$ 
such that $p_\alpha \leq p\cup \{ \langle \beta _\alpha ,\beta _\alpha 
\rangle \} ,$ $p_\alpha \Vert \dot q_\alpha$ if $\beta _\alpha \in S.$ Define
$f:S\cap \{ \beta _\alpha :\alpha <\omega _1\} \to \omega _1$ by 
$f(\beta _\alpha )=g(p_\alpha \cap M_\alpha ).$
$f$ is defined on a stationary set and can be easily seen to be 
regressive, therefore we can find a stationary set 
$T\subset \omega _1$ and $q$ such that 
$\alpha \in T$ implies $p_\alpha \cap M_\alpha=q.$ Now 
similarly to the proof of the $\Delta$-system lemma one 
can find $U\subset T$ of cardinality $\aleph _1$ such 
that $\alpha _1,\alpha _2\in U$ implies $p_{\alpha _1}$ 
is comparable with $p_{\alpha _2}$ and thus $\{ \dot q_\alpha /
p_\alpha :\alpha \in U\} \subset Q$ is an antichain giving 
contradiction with assumed c.c.c. of $Q.$
\enddemo

\proclaim
{Lemma 4} Let $S\subset \omega _1$ be stationary and $J=NS_{\omega _1}+
(\omega _1\setminus S).$ \roster
\item If $J$ is precipitous then $P_S\Vdash$``$NS_{\omega _1}$ is precipitous".
\item If $J$ is presaturated then $P_S\Vdash$``$NS_{\omega _1}$ is 
presaturated" iff $\Cal P(\omega _1)/I\Vdash$``$jS\subset \omega _2^V=
\omega _1^{V[G]}$ is stationary", where $j:V\to M$ is the canonical 
generic ultrapower.
\endroster
\endproclaim

{\it Remark.} The situation in (2) is parallel to that of \cite {BT}.
 Notice that $P_S\Vdash$``$NS_{\omega _1}$ is not $\omega _2-$saturated". 
( Choose $\langle g_\gamma :\gamma <\omega _2\rangle ,$ a family of 
almost disjoint functions, $g_\gamma :\omega _1\to S,$ $g_\gamma 
(\alpha )>\alpha ,$ all $\alpha <\omega _1, $ $\gamma <\omega _2.$ 
If $G:\omega _1\to S$ is the $P_S$-generic function then $\langle 
S_\gamma :\gamma <\omega _2\rangle ,$ $S_\gamma =\{ \alpha <\omega _1:
 G(\alpha )=g_\gamma (\alpha )\}$ is a long antichain of stationary sets 
in $V[G].)$ As far as the second condition in (2) is concerned, let us say 
that  a pair $Q,j$ is stationarily correct if $Q\Vdash$``$j:V\to M$ 
is elementary, $\kappa =crit(j)$ and $M\models$``$T\subset j\kappa$ 
is stationary" iff $T$ is stationary". We have proved that the 
nonstationary tower ultrapower 
as described in \cite {W} is stationarily correct as well as the 
$NS_{\omega _1}$-generic ultrapower under $MA^+(\omega _1$-closed). 
It is consistent w.r.t. suitable hypotheses that $NS_{\omega _1}$ is 
c.c.c. destructibly $\omega _2$-saturated and stationarily correct or 
that there is $J\subset \Cal P(\omega _1),$ a normal $\omega _2$-saturated 
ideal which is not stationarily correct. It seems however that it is an 
open problem whether $NS_{\omega _1}$ plus a single set can be presaturated 
and not stationarily correct. Thus the second condition in (2) 
is possibly empty.

\demo
{Proof} (1) follows from the following two claims:

\proclaim
{Fact 1} $P$ preserves maximal antichains of stationary subsets of $S.$
\endproclaim

\proclaim
{Fact 2} Let $p\in P,$ $p\Vdash$``$\dot f:\dot T_0\to Ord,$ $\dot T_0
\subset \omega _1$ stationary. Then there are $q<p,\dot T_1,g:\omega _1
\to Ord$ such that $q\Vdash$``$\dot T_1\subset \dot T_0$ is stationary 
and $\forall \alpha \in \dot T$ $\dot f(\alpha )=\check g(\alpha )".$
\endproclaim

Let us fix an enumeration $E:P_S\to \omega _1$ and go on to prove the 
above facts. In the case of Fact 1, let $\langle T_i:i\in I\rangle$ be 
a maximal antichains of stationary subsets of $S,$ $p\in P_S,$ $p\Vdash$
``$\dot U\subset \omega _1$ is stationary". Set $\bar U=\{ \alpha \in 
S:\exists q\leq p\cup \{ \langle \alpha ,\alpha \rangle \}$ $q\Vdash 
\alpha \in \dot U\} $ and choose $q_\alpha$ witnessing $\alpha \in \bar U$ 
for $\alpha$ in $\bar U.$ $\bar U$ is stationary. Define $F:\bar U\to 
\omega _1$ by $F(\alpha )=E(q_\alpha \restriction \alpha ).$ $F$ is 
regressive on a stationary set and so we can find $T\subset \bar U$ 
stationary such that $F^{\prime \prime }T=\{ \beta \} ,$ some $\beta 
<\omega _1 .$ Let $q=F^{-1}(\beta )$ and choose $i\in I$ such that 
$T\cap T_i$ is stationary. Then $p>q\Vdash$``$\dot U\cap \check T_i$ 
is stationary": let $r<q,$ $r\Vdash$``$\dot C\subset \omega _1$ is a 
club". Choose $M\prec H_\theta$ countable for some large regular $\theta$ 
such that $r,\dot C\in M$ and $\alpha =M\cap \omega _1\in T.$ Then 
$q_\alpha$ and $r$ are compatible and their common lower bound forces
 $\alpha$ into $\dot U\cap \check T_i\cap \dot C$ (notice $q_\alpha$ 
is a master condition for $M).$

The proof of Fact 2 follows a similar pattern. Let $p,T_0$ be as in the 
Fact 2. Set $\bar T_0=\{ \alpha \in S:\exists  q\leq p\cup \{ \langle 
\alpha ,\alpha \rangle \}$ $q\Vdash \alpha \in \dot T_0\} .$ For $\alpha 
\in \bar T_0$ choose $q_\alpha \leq p\cup \{ \langle \alpha ,\alpha \rangle 
\}$ $q_\alpha \Vdash \alpha \in \dot T_0,$ $q_\alpha$ decides $\dot 
f(\alpha ).$ Let $F(\alpha )=E(q_\alpha \restriction \alpha ).$ $F$ 
is regressive on a stationary set and we can find $\beta ,\bar T_1\subset
 S$ stationary, $F^{\prime \prime }\bar T_1=\{ \beta \} .$ Define $\dot T_1=
\{ \alpha \in \bar T_1: q_\alpha \in G\} ,$ where $G$ is the generic filter 
and $g:T_1\to Ord,$ $g(\alpha )=$ the unique $\xi$ such that $q_\alpha 
\Vdash \dot f(\alpha )=\check \xi .$ Then as above $p>F^{-1}
(\beta )\Vdash$``$\dot T_1$ is stationary and $\dot f\restriction 
\dot T_1=\check g\restriction \dot T_1".$

It is left to the reader to show now that (1) holds. To prove (2) 
we first observe

\proclaim
{Fact 3} $\Cal P(\omega _1)/J*jP_S\restriction \langle \omega _1,
\omega _1 \rangle$ is isomorphic to $P_S*\dot \Cal P(\omega _1)/
NS_{\omega _1}.$
\endproclaim

To see this, let $G\subset \Cal P(\omega _1)/I$ be $V$-generic and 
$H\subset jP_S\restriction \langle \omega _1,\omega _1 \rangle$ 
be $V[G]$-generic. Again we confuse $H$ with $\bigcup H:\omega _2^V\to jS.$ 
Set $H^\prime =H\restriction \omega _1.$ Thus $H^\prime$ can be regarded 
as $V$-generic object for $P_S.$ The standard techniques give an extension 
of $j:V\to M,$ $j\in V[G],$ to $\hat j:V[H^\prime ]\to M[H]$ in $V[G][H].$ 
We set $G^\prime =\{ T\in \Cal P(\omega _1)\cap V[H^\prime ]:\omega _1 \in 
\hat jT\}$ and claim that $G^\prime \subset \Cal P(\omega _1)/NS_{\omega _1}$ 
is $V[H^\prime ]$-generic and moreover $V[G][H]=V[H^\prime ][G^\prime ].$ 
To this end, fix $f,T,p,\langle \dot T_i:i\in I\rangle$ such that 
$T\subset S$ is stationary, $T\Vdash$``$\langle \omega _1,\omega _1\rangle 
\in \check p\in jP_S,$ $\check p=[\check f]",$ $f:T\to P_S,$ 
$p\restriction \omega _1\Vdash _{P_S}$``$\langle \dot T_i:i\in I\rangle$ 
is a maximal antichain in $\Cal P(\omega _1)/NS_{\omega _1}".$ For $i\in 
I$ set $\bar T_i=\{ \alpha \in S:\exists q\in P_S$ $q\leq f(\alpha ),$
 $q\Vdash \alpha \in \dot T_i\} .$ Since the $\dot T_i$'s are forced to form 
a maximal antichain, there is an $i\in I$ such that $\bar T_i\cap T$ 
is stationary. For each $\alpha \in \bar T_i\cap T,$ choose $q_\alpha 
\leq f(\alpha ),$ $q_\alpha \Vdash \alpha \in \dot T_i\} .$ Define 
$F(\alpha )=E(q_\alpha \restriction \alpha ).$ $F$ is regressive on 
a stationary set and we can find $U\subset \bar T_i\cap T,\beta <\omega _1$ 
such that $U$ is stationary and $F^{\prime \prime }U=\{ \beta \} .$
Let $g:U\to P_S$ be defined by $g(\alpha )=q_\alpha .$ Then in $\Cal
 P(\omega _1)/I*jP_S\restriction \langle \omega _1,\omega _1 \rangle$ $U,[g]
\leq T,p$ and $U,[g]\Vdash \omega _1 \in \hat j(\dot T_i/H^\prime )$ 
and therefore $G^\prime$ meets the antichain given by $\langle \dot 
T_i:i\in I\rangle .$ This proves the genericity. To reconstruct $G,H$
 from $G^\prime ,H^\prime ,$ notice that $G=G^\prime \cap V.$ If $\hat j$
 is the generic ultrapower of $V[H^\prime ]$ by $G^\prime ,$ it is
immediate that $H=\hat jH^\prime .$

(2) now follows: if $jS$ is forced to be stationary, then $jP_S$ will 
be a forcing in $V[G]$ which does not collapse $\omega _1^{V[G]}=
\omega _2^V.$ If on the other hand $jS$ can be nonstationary, let us 
say $T\Vdash$``$jS$ is nonstationary" then it is easy to find two 
disjoint closed unbounded subsets of $\omega _2^V$ in $V[G][H]=
V[H^\prime ][G^\prime ]$ if $T\in G\subset G^\prime$ and so 
$\omega _2$ was collapsed.

\enddemo

From now on, let $Q_{\omega _1}$ denote the forcing adding $\omega _1$
 Cohen reals. From our previous work, any real added by $P_S$ is in 
some Cohen extension of the ground model. It is also not very hard 
to see that $Q_{\omega _1}$ regularly embeds into $P_S.$ (See Corollary 5 
for a rather complicated example how to do this.) It is natural to ask 
whether such embedding can reap all the real numbers of $V^{P_S},$ i.e. 
if we can have $Q_{\omega _1}\lessdot P_S$ as a regular subalgebra so 
that $P_S\Vdash$``$\omega ^\omega \cap V[G]=\omega ^\omega \cap 
V[G\cap Q_{\omega _1}]".$

\proclaim
{Lemma 5} Let $H\subset Q_{\omega _1}$ be generic. In $V[H]$ (actually 
in $V[\omega ^\omega ]^{V[H]})$ there is an $\omega$-distributive 
$S$-proper  forcing $T$ such that $T\Vdash$``there is $G\subset P_S$ 
$V$-generic such that 
$\omega ^\omega \cap V[G]=\omega ^\omega \cap V[H]".$
\endproclaim

\demo
{Proof} Work in $V[H]$ and define $T=\{ g:\exists \alpha <\omega _1,$
 $\{ \langle \alpha ,\alpha \rangle \} \in P_S,$ $g\subset P_{S,\alpha }$
 is generic over $V\}$ ordered by reverse inclusion. Certainly all $g\in T$ 
are hereditarily countable, thus coded by reals and $T\in V[\omega ^\omega ].$
 Choose $g\in T,$ $\langle D_i:i<\omega \rangle$ a sequence of open 
dense subsets of $T,$ and $M\prec H_\theta ,$ $S,g,\langle D_i:i<\omega 
\rangle ,H\in M,$ $M\cap \omega _1\in S.$ $\{ T\cap M\}\cup \{ D\cap 
T\cap M:D\in M$ open dense $\}$ is a countable collection of 
hereditarily countable objects and as such belongs to some 
$V[H\restriction \alpha ],$ $\alpha <\omega _1.$ In $V[H\restriction 
\alpha ],$ $T\cap M\restriction g$ is isomorphic to adding one Cohen 
real. Let us regard $H(\alpha )$ as a subset of it. Then it is easy to
 show that $h=\bigcup H(\alpha )$ is a $V$-generic subset of $P_{S,M\cap 
\omega _1}$ and a strongly $M$-generic condition under $g,$ in particular
 $h\in \bigcup _{i<\omega }D_i.$

Due to the local genericity condition in Corollary 4, if $K\subset T$ 
is generic, $G=\bigcup K$ is a $V$-generic subset of $P_S.$ The last thing 
to check is that $\omega ^\omega \cap V[G]=\omega ^\omega \cap V[H].$ To 
this aim, for $\alpha <\omega _1$ define $D_\alpha =\{ g\in T: 
H\restriction \alpha \in V[g]\} .$ The following Subclaim will 
complete the proof.

\proclaim
{Subclaim} Each $D_\alpha$ is a dense subset of $T.$
\endproclaim

\enddemo

Thus the reals coming from $P_S$ look exactly the same as the reals 
coming from $Q_{\omega _1}.$

\proclaim
{Lemma 6} Cons(ZFC+$\kappa$ Mahlo) implies Cons(ZFC+$\exists S\subset
 \omega _1$ stationary costationary and there is an embedding 
$Q_{\omega _1}\lessdot P_S$ reaping all the reals of $V^{P_S}.$
\endproclaim
 
\demo
{Proof} Fix a Mahlo cardinal $\kappa$ and set $S=\{ \alpha <\kappa :
\alpha$ inaccessible $\}.$ $Coll(\omega ,<\kappa )$ is homogeneous 
and so for every finite function $p$ from $\kappa$ to $\kappa ,$ 
either $Coll(\omega ,<\kappa)\Vdash$``$\exists \alpha<\kappa$ $\exists f:
\alpha \to \kappa$ increasing continuous with $f^{\prime \prime }\alpha 
\subset S,$ $p\subset f"$ or $Coll(\omega ,<\kappa)\Vdash$``$\lnot 
\exists \alpha<\kappa$ $\exists f:\alpha \to \kappa$ increasing continuous
 with $f^{\prime \prime }\alpha \subset S,$ $p\subset f".$ (Notice that due 
to the $\kappa$-c.c. $Coll(\omega ,<\kappa )$ preserves stationarity of 
$S.)$ Therefore we can define $P=\{ p:p$ is a finite function and 
$Coll(\omega ,<\kappa)\Vdash$``$\exists \alpha<\kappa$ $\exists f:
\alpha \to \kappa$ increasing continuous with $f^{\prime \prime }\
alpha \subset S,$ $p\subset f"$ ordered by inclusion and be sure to
 get $Coll(\omega ,<\kappa )\Vdash$``$\check P=\dot P_{\check S}.$

\proclaim
{Claim 1} $P\Vdash$``$\check \kappa =\dot \aleph _1,G$ (confused with 
$\bigcup G):\kappa \to \check S$ increasing continuous."
\endproclaim

\proclaim
{Claim 2} $P\Vdash$``$\exists H\subset Coll(\omega ,<\kappa )$ generic 
over $V,$ $\omega ^\omega \cap V[H]=\omega ^\omega \cap V[G]"$
\endproclaim

\demo
{Proof} Fix $G\subset P$ generic and work in $V[G].$ Notice that as in the 
case of Lemma 1, (due to the easy factorization of $P)$ $r\in 
\omega ^\omega \cap V[G]$ iff $r\in \omega ^\omega \cap V[G\restriction 
\alpha ],$ some $\alpha <\kappa .$ Consider the following poset $Y=\{
 h: h\subset Coll(\omega  ,<\alpha )$ generic over $V$ for some $\alpha 
\in rng(G)\} $ ordered by reverse incusion. $Y$ is $\omega$-closed by 
the closure of $rng(G):$ assume $h_0>h_1>\dots >h_i>\dots ,$ $i<\omega ,$ 
is a decreasing sequence of elements in $Y,$ $h_i\subset Coll(\omega ,
<\alpha _i),$ some $\alpha _i\in rng(G).$ Then $\alpha =sup_{i<\omega }
\alpha _i\in rng(G),$ $\alpha$ is $V$-inaccessible and $h=
\bigcup _{i<\omega }h_i\subset Coll(\omega ,<\alpha )$ is generic over 
$V,$ since if $A\subset Coll(\omega ,<\alpha )$ is a maximal antichain 
in $V,$ we have $|A|<\alpha$ (in $V)$ and thus for some $i<\omega $ 
$A\subset Coll(\omega ,<\alpha _i )$ and $A$ is met by a condition in
 $h_i.$ For $\alpha <\kappa$ define $D_\alpha =\{ h\in Y:G\restriction 
\alpha \in V[h]\} .$ The following subclaim will finish the proof of the 
Claim 2 since $\kappa =\aleph _1,$ $Y$ is $\omega $-closed and any real 
in $V[G]$ is coded by an initial segment of $G.$

\proclaim
{Subclaim} Each $D_\alpha$ is dense in $Y.$
\endproclaim
\enddemo

Now we can finish the proof of the Lemma. Fix $\dot H,$ a $P$-name
 for a generic subset of $Coll(\omega ,<\kappa )$ as in Claim 2. 
Fix $K\subset Coll(\omega ,<\kappa )$ generic over $V.$ We claim 
that $V[K]$ is a model of the wanted theory with our $S.$ To prove 
it, choose $G\subset P_S=P$ generic over $V[K].$ By a mutual 
genericity argument, $H=\dot H/G\subset Coll(\omega ,<\kappa )$ is
 generic over $V[K].$ In $V[K],$ however, $\kappa =\aleph _1$ and 
so $Coll(\omega ,<\kappa )$ is isomorphic to $Q_{\omega _1}.$ We 
view $H$ as a subset of $Q_{\omega _1}$ (transferred by some 
isomorphism of $Coll(\omega ,<\kappa )$ and $Q_{\omega _1}$ in $V[K]).$
 The only thing left to check is that $\omega ^\omega \cap V[K][G]=
\omega ^\omega \cap V[K][H].$ For that we use Corollary 3 and the 
most significant property of $H,$ that $\{ G\restriction \alpha :
\alpha <\kappa\} \subset V[H].$ 
\enddemo

\proclaim
{Lemma 7} If $C\subset \omega _1$ is a club then $P_{S\cap C}\lessdot 
P_S.$ In fact, $P_S\Vdash$``if $\dot D\subset \omega _1$ is the generic 
club then $\check C\cap \dot D$ is $P_{S\cap C}$-generic club".
\endproclaim

\demo
{Proof} By the local genericity criterion in Corollary 4 it is enough 
to prove the following claim:
\proclaim
{Claim 3} If $\alpha <\omega _1$ is idecomposable, $T\subset S\subset 
\alpha ,$ $T$ clunbounded in $S$ and $P_{T,\alpha },P_{S,\alpha }$ are
 nonempty posets then $P_{S,\alpha }\Vdash$``if $\dot D\subset \alpha$
 is the generic club then $\dot D\cap T$ is a $P_{T,\alpha }$-generic club".
\endproclaim
\demo
{Proof of the Claim} We first give two subclaims, then prove the Claim 
from them and complete the proof of the Lemma by proving the two subclaims.
\proclaim
{Subclaim} If $\gamma$ is indecomposable, $T\subset S,$ where $T$ is 
clunbounded in $S,$ which is countable and $P_{T,\gamma },P_{S,\gamma }
\neq 0$ then
$P_{S,\gamma }\Vdash$``$o.t.\dot D\cap \check T=\gamma",$ where $\dot D$
 is the generic club through $S.$
\endproclaim
\proclaim
{Subclaim} If $I,J\subset \omega _1$ are countable intervals of 
ordinals, $o.t.J\leq o.t.I$ are both indecomposable, $T\subset S,$
 where $T$ is clunbounded in $S,$ which is countable and $P_{T,J },
P_{S,I }\neq 0$ then
for any $t\in P_{T,J}$ there is $s\in P_{S,I}$ such that $s\Vdash$
``$t$ is a subset of the increasing enumeration of $\check T\cap 
\dot D$ starting with $min(I)",$
 $\dot D$ the generic club.

\endproclaim

Now we can proceed to prove the Claim. For technical reasons we pretend 
that $\alpha \in S\cap T$ and any $p\in P_{S,\alpha }$ contains $\langle 
\alpha ,\alpha \rangle$ (accordingly $\alpha \in \dot D$ then). 
Choose $p_0\in P_{S,\alpha }$ arbitrary. We find $p <p_0$ and $q\in 
P_{T,\alpha }$ such that for any $q^\prime \leq q$ there is $p^\prime
 <p$ such that $p^\prime \Vdash _{P_{S,\alpha }}$``$q^\prime \subset$ 
the enumeration of $\dot D\cap T",$ proving the Claim.  We build 
$\alpha =\alpha _0>\alpha _1 >\dots ,$ $p_0=p_0>p_1>\dots ,$
 $\alpha =\gamma _0>\gamma _1>\dots$ so that
\roster
\item $\alpha _i\in dom(p_i),$
\item $p_i\Vdash _{P_{S,\alpha }}$``$p_i(\alpha _i)\in T$ is the
 $\gamma _i^{th}$ element of $\check T\cap \dot D"$ where $\dot D$
 is the generic club $\subset \alpha ,$
\item $dom(p_{i+1}\setminus p_i)\subset \alpha _i,$
\item $o.t.\alpha _i\setminus \alpha _{i+1}$ is indecomposable,
\item $o.t.\gamma _i\setminus \gamma _{i+1}$ is indecomposable,
\item $\gamma _i$ limit implies $dom(p_{i+1})\cap (\alpha _{i+1},\alpha _i)=0,$
\item $\gamma _i$ successor implies $\gamma _{i+1}$ is the 
predecessor of $\gamma _i.$
\endroster
This is easily done and must end at some $n<\omega$ since the $\alpha _i$'s
 form a descending sequence of ordinals. Set $q=\{ \langle \gamma _i,
p_i(\alpha _i)\rangle :i<n\}$ and $p=p_n.$ We claim that $p,q$ are what
 we are looking for. First, $q\in P_{T,\alpha }.$ Choose $M\prec 
H_\theta$ countable containing everything relevant and $p\in 
g\subset P_{S,\alpha }$ generic over $M.$ Then by elementary 
absoluteness considerations $T\cap rng(\bigcup g)$ is a club subset of
 $T$ of ordertype $\alpha$ such that $q$ is a subset of its enumeration.
 Second, choose $q^\prime <q$ in $P_{T,\alpha }.$ Let us assume for
 simplicity that $dom(q^\prime \setminus q)\subset (\gamma _{i+1},
\gamma _i)$ for some $i.$
Then $\gamma _i$ is limit by (7). Now we use the second subclaim with 
$I=\alpha _i\setminus \alpha _{i+1},$ $J=\gamma _i\setminus 
\gamma _{i+1},$ $S\cap (p_{i+1}(\alpha _i),p_i(\alpha _i))$ 
in the place of $S$ and $T\cap (p_{i+1}(\alpha _i),p_i(\alpha _i))$
 in the place of $T$ on $t=q^\prime \setminus q.$ The resulting $s$ 
is easily seen to be such that $P_{S,\alpha }\ni p\cup s=p^\prime
\Vdash _{P_{S,\alpha }}$``$q^\prime$ is a subset of the increasing 
enumeration of $\dot D\cap \check T".$

To prove the first subclaim, let $\gamma <\omega _1$ be the least 
indecomposable such that there are $T\subset S$ violating the statement. 
We distinguish two cases:\roster
\item $\gamma$ is a limit of indecomposables. Let $P_{S,\gamma }
\ni p\Vdash$``$o.t.\dot D\cap \check T<\beta <\gamma "$ for some
 indecomposable $\beta >max(dom(p)).$ Let $\xi \in T,$ $\xi >max(rng(p))$
 be such that $P_{T\cap \xi ,\beta }\neq 0.$ By indecomposability of
 $\beta$ and $\gamma$ $p\cup \{ \langle \beta ,\xi \rangle \} \in 
P_{S,\gamma }.$ By minimality of $\gamma ,$ $P_{S\cap \xi ,
\beta }\Vdash$``$o.t.\check T\cap \xi \cap D=\beta ",$ 
therefore  
$p\cup \{ \langle \beta ,\xi \rangle \} \Vdash$``$o.t.
\check T\cap \dot D\geq \beta ",$ a contradiction.
\item $\gamma =\omega \beta$ for some $\beta <\gamma$ 
indecomposable. Fix $p\in P_{S,\gamma }.$ We get $q=
p\cup \{ \langle \delta ,\xi \rangle \} <p$ such that $q\Vdash$
``$o.t.\check T\cap \dot D\cap (\xi \setminus max(rng(p)))=\beta 
".$
The contradiction then follows by a simple genericity argument. 
To get our $q,$ we set $\delta =max(dom(p))+\beta$ and choose 
$\xi \in T,$ $\xi >max(rng(p))$ such that $P_{T\cap \xi ,\beta }\neq 0.$ 
By minimality of $\gamma$ and indecomposability of $\gamma ,\beta ,$ it 
follows that $q$ works.
\endroster

The second subclaim is in fact a corollary to the first one. It is 
certainly enough to prove it for $J=\gamma \leq \alpha =I.$ For 
simplicity we assume that $t=\langle \beta ,\xi \rangle ,$ $\xi 
\in T.$ Find $\beta =\beta _0>\beta _1>\dots >\beta _{m-1}>0=\beta _m$ 
so that $\forall i<m$ $o.t.(\beta _i\setminus \beta _{i+1})$ is 
indecomposable and choose $s\in P_{T,\gamma },$ $dom(s)=
\{ \beta _0\dots \beta _{m-1}\} ,$ $s(\beta )=\xi .$ Then $s$ is a 
member of $P_{S,\alpha }$ as well and by the first claim it has the 
required property.

\enddemo
\enddemo

Let us evaluate the factor forcing $P_S/P_{S\cap C}.$ For definiteness, 
assume that $C\subset \omega _1$ is such that $|S\setminus C|=\aleph _1.$ 
Let us fix $H\subset P_{S\cap C}$ $V$-generic and in $V[H],$ choose a
 continuous increasing sequence $\langle M_\alpha :\alpha <\omega _1\rangle$ 
of countable submodels of some $H_\theta$ with $C,S,H\in M_0.$ Define
 $\beta _0=0$ and for $0<\alpha <\omega _1$ let $\beta _\alpha =M_\alpha
 \cap \omega _1.$ Define a forcing $Q=$ the finite support product of 
$P_{S\cap I_\alpha }/P_{S\cap C\cap I_\alpha}$ for $\alpha <\omega _1$
 where $I_\alpha =[\beta _\alpha ,\beta _{\alpha +1}),$ the $V$-generic 
subset of $P_{S\cap C\cap I_\alpha}$ is just 
$H\cap P_{S\cap C\cap I_\alpha}$ and the embedding
 $P_{S\cap C\cap I_\alpha}\lessdot P_{S\cap I_\alpha}$ is the one described 
in Claim 3 (modulo an ordinal shift) . Then it is not difficult to see that 
$Q$ is isomorphic to $Q_{\omega _1}$ in $V[H],$ since the forcings standing
 in the finite support product are nontrivial and $\aleph _0$-dense. One
 can easily prove that a $V[H]$-generic $K\subset Q$ together with $H$ gives
 a $V$-generic $G\subset P_S$ such that $V[G]=V[H][K].$

\proclaim
{Corollary 5} $P_S=P_S\times Q_{\omega _1}.$
\endproclaim
Since $C\subset \omega _1$ is as above, $P_S=P_{S\cap C}*Q_{\omega _1}
=P_{S\cap C}\times Q_{\omega _1}=P_{S\cap C}\times Q_{\omega _1}\times 
Q_{\omega _1}=P_S\times Q_{\omega _1}.$

\proclaim
{Corollary 6} If $S=T$ modulo $NS_{\omega _1}$ then $P_S=P_T$ (again, 
as Boolean algebras).
\endproclaim
To see this, fix $S,T$ as in the Corollary and choose $C\subset \omega _1$ 
club such that $|S\setminus C|=|T\setminus C|=\aleph _1$ and $S\cap C=T\cap 
C.$ Then $P_S=P_{S\cap C}\times Q_{\omega _1}=P_{T\cap C}\times 
Q_{\omega _1}=P_T$ from the remarks preceding Corollary 5.

\proclaim
{Corollary 7} Cons($\kappa$ Mahlo) implies Cons($Q_{\omega _1}$ embeds 
into $P_{\omega _1}$ reaping all the reals).
\endproclaim

The proof of the corollary is left to the reader. The model is $V[K][L],$
where $K\subset Coll(\omega ,<\kappa)$ is generic as in Lemma 6 and $K$ 
is a generic club through $S=\{ \alpha <\kappa :\alpha$ inaccessible in 
$V\}$ using {\it countable} conditions. The key to the proof is to notice 
that $V[K][L]\models P_{\omega _1}=P_S$ as Boolean algebras; the rest 
carries over from Lemma 6.

\enddocument